\documentclass{amsart}
\usepackage{amsmath,amssymb,enumerate,graphicx,cases}

\renewcommand{\P}{\mathbb{P}}
\renewcommand{\L}{\mathcal{L}}

\providecommand{\E}{\mathbb{E}}
\providecommand{\F}{\mathcal{F}}
\providecommand{\A}{\mathcal{A}}

\providecommand{\R}{\mathbb{R}}
\providecommand{\C}{\mathcal{C}}

\providecommand{\N}{\mathbb{N}}
\providecommand{\Z}{\mathbb{Z}}

\newtheorem{theo}{Theorem}[section]
\newtheorem{coro}[theo]{Corollary}
\newtheorem{lem}[theo]{Lemma}

\newtheorem{prop}[theo]{Proposition}
\newtheorem{rem}[theo]{Remark}

\newtheorem{defi}[theo]{Definition}

\begin{document}

\title{\textbf{ Exit time for anchored expansion}}
\author{Thierry Delmotte,  Cl\'ement Rau}
\date{}

\begin{abstract}
Let $(X_n)_{n\geq 0}$ be a reversible random walk 
on a graph $G$ satisfying an anchored isoperimetric
inequality. We give upper bounds for exit time
(and occupation time in transient case) by X
of any set which contains the root. As an application,
we consider random environments of $\Z^d$.
 
\end{abstract}

\maketitle
\tableofcontents
\section{Introduction}
\begin{rem} A shorter version of this paper is proposed to Annales de la Facult{\'e} des Sciences de
Toulouse without the proof of isoperimetry for
random environment (Proposition 5.2), which is improved since the
first version on Arxiv. This proof will be proposed in
another article, in the second time. 
\end{rem}
Among many connections linking the geometry of a graph
and the behaviour of its simple random walk, one
important tool is isoperimetry.
Already present in the
celebrated work of Nash \cite{nash}, the idea was made
clear since the seminal work of Varopoulos \cite{varo}.
A discrete version for graphs is in \cite{cosa}.
The isoperimetry conditions are various, geometric or
functional. For instance the inequality above yields a
$L^2$ type of Faber-Krahn inequality proposed by
Grigor'yan in \cite{grigFK}, and Coulhon has shown
in \cite{coulhon} that it gives an upper bound of the
iterated transition probabilities of the random walk.

The problem
with uniform isoperimetric inequality is its unstability
under random perturbations like percolation.
If one studies the ''ant in the labyrinth'' of
de Gennes \cite{degennes}, one needs a weaker version of
isoperimetry which can be robust, as introduced in the two
last decades by Thomassen in \cite{thomassen} and next by
Benjamini, Lyons and Schramm in \cite{benj}. It is called
{\it{ anchored or rooted isoperimetric inequality}}. Here
is the definition.
%graphs that we consider here   have not necessarly  bounded degrees and  
For a graph $G$, we denote $V(G)$ the set of vertices and
$E(G)$ the set of edges.
\begin{defi}Let $\F$ a positive increasing function defined on $\R_+$. Let  $G$ a graph and $ o\in G$.  We say that 
$G$ satisfies an anchored (or rooted) $\F$-isoperimetric inequality at o if  there exists a constant $C_{\rm IS} >0$ such that for any connected set $A$ which contains o we have:
\begin{equation} \label{IS}
\frac{ |\partial A|}{ \F(|A|)}\geq  C_{\rm IS}.
\end{equation}
$\partial A$ is equal to the set $\{(x,y) \in E(G); \ x\in A \ and \ y\notin A\}$ and $|B| $ stands for the cardinal of $B$.
\end{defi}
When $\F(x)=x^{1-1/d}$, we will say that $G$ satisfies a $d-$dimensional isoperimetric inequality.
When $\F=id$ and $G$ has bounded degree there is an equivalent
version of this definition which reads as follows: \\
  {\it{ $G$ satisfies a strong anchored (or rooted) isoperimetric inequality  if 
  $$\lim_{n\rightarrow \infty} \inf\{ \frac{|\partial S|}{|S|}; \ S \ connected, \ v\in S \ and \ |S|\geq n \}:=i(G)$$ is strictly positif.}}

This definition does not depend on the choice of the fixed vertex whereas in the previous definition, the constant $ C_{\rm IS}$ depends on the point o.

Our object here is to examine what anchored isoperimetric
inequality implies for random walk. Our hope is that it could
be useful for instance in the still very open problem of
spectral dimension which could not be $4/3$ in low dimensions
as stated by the Alexander-Orbach conjecture \cite{ao}.
See the lecture of Barlow \cite{barlow} for an introduction
to these questions.

\subsection{What we know for anchored expansion.}$\ $ \\
The first result known for rooted $\F$-isoperimetric inequality is due to Thomassen. In \cite{thomassen}, it is proved that a the simple random walk on a graph $G$ is transient if $G$ satisfies a rooted $\F$-isoperimetric inequality such that $\sum_k \F(k)^{-2} <\infty$ . The main step of the proof is  to extract a subdivision of the dyadic tree from the initial graph. Then, thanks to hypothesis, it is possible to construct a finite flow  on the tree, which proves that the  tree is transient.

It was long afterwards that other results did appear for anchored expansion.
In 2000 Virag has studied the case of strong anchored isoperimetric inequality.
In \cite{virag}, it is proved that strong anchored isoperimetric  inequality on graphs with  bounded geometry, implies a positive lim inf speed. Moreover Virag proves  that in this case, transitions probability at time $n$ of the random walk are bounded by $e^{-n^{1/3}}$.
 
Later, still when $\F=id$, Chen and Peres have proved that if $G$ satisfies a strong anchored isoperimetric inequality then so does every infinite cluster of independant percolation with parameter $p$ sufficiently closed to $1$. Next, they have shown that strong anchored expansion is preserved  under a random stretch if, and only if, the stretching law has an exponential tail. They also proved that for a supercritical Galton Watson tree $\mathbb{T}$ given nonextinction, we have $i(\mathbb{T})>0$ a.s. \\

  \subsection{What we don't know. Open questions}$\ $ \\
  There is an important collection of conjectures relating to anchored expansion. Here is some of them: \\ 
Question 1:  does the sub tree of Thomassen satisfy an anchored Isoperimetric inequality ? \\
Question 2: does a general anchored isoperimetric inequality imply an upper bound of  $p_n(x,y)$ ? \\
Question 3: does anchored isoperimetric inequality is a good tool to prove an invariance principle in random environment of $\Z^d$ ? \\

\subsection{Continuous space setting.} 

The paper is written in the discrete space setting of graphs. The
reason is that anchored isoperimetric inequality is a natural tool
in random media and is therefore more associated with this setting.
In fact the continuous setting (of Riemannian manifolds for instance)
works as well, and may be, the proofs are far more readable.
As both an introduction to our technique and an illustration of what
the continuous setting results would look like, we begin with a key
result written in this setting. Details, especially from potential
theory, will only appear later in the paper for graphs.

Let $M$ be a Riemannian manifold with an anchored isoperimetric
inequality at root $o$, that is (\ref{IS}) for finite volume smooth connected
domains $A$ containing $o$. Precisely, $|A|=m(A)$ for the Riemannian
volume element $m$ and $|\partial A|=\mu(\partial A)$ for the
Riemannian volume element $\mu$ on the smooth submanifold $\partial A$. \\
Now let fix some $A$ and consider the Brownian motion on $M$ starting
at $o$ and killed when hitting $\partial A$ at time $\tau_A$.
We denote $p_t^A$ its submarkovian
kernel, $A_s$ the level sets of Green function and $u(s)$ their measures.
$$A_s=\left\{x\in A, G^A(x)=\int_0^\infty p_t^A(x)\ {\rm d}t \geq s\right\},
\qquad u(s)=m(A_s).$$
Thanks to harmonic properties of $G^A$, these level sets are connected
and contain the root. Thus, they will also satisfy (\ref{IS}).
In the following we use $\mu$ for any $s$ and also $\nu$ denoting the
inward unit normal vector field on $\partial A_s$. The inward direction
is chosen to have $G^A$ increasing.

\begin{theo} \label{thcont}
The anchored isoperimetric inequality yields a differential inequation
$$u'(s)\leq -\Big(C_{\rm IS}\F(u(s))\Big)^2.$$
This naturally leads to upper estimates of $u(s)$ and
$\E(\tau_A)=\int_0^\infty u(s)\ {\rm d}s$. \\
For instance if $\F(u)=u^{1-1/d}$, $\E(\tau_A)\leq C m(A)^{2/d}$,
\newline and if $\F(u)=u$, $\E(\tau_A)\leq C \ln m(A)$.
\end{theo}

\begin{proof}
Schwarz inequality
$$\Big(C_{\rm IS}\F(u(s))\Big)^2\leq\mu(\partial A_s)^2
=\left(\int_{\partial A_s}{\rm d}\mu\right)^2
\leq \int_{\partial A_s}\frac{\partial G^A}{\partial\nu}{\rm d}\mu
\int_{\partial A_s}\frac{{\rm d}\mu}{\partial G^A/\partial\nu}$$
involves the flow
$$\int_{\partial A_s}\frac{\partial G^A}{\partial\nu}{\rm d}\mu=1$$
and the derivative of $u$ since whith the co-area formula,
$$u(s)=\int_{G^A\geq s}{\rm d}m=\int_s^\infty\left(\int_{G^A=t}
\frac{{\rm d}\mu}{\partial G^A/\partial\nu}\right){\rm d}t.$$
This yields the differential inequation.

For $\F(u)=u^{1-1/d}$, computations may be avoided if we
compare with the case when $A$ is a ball of radius $R$ in
$\mathbb{R}^d$. Then $\partial G^A/\partial\nu$ is constant,
all inequalities are equalities and the result should be
that $\E(\tau_A)$ is like $R^2$.
\end{proof}

Application of this Schwarz inequality is already
apparent in \cite{grig}, \cite{hesch} or \cite{lms}
to establish a recurrence criterion or estimate resistance.

 \subsection{Results of the paper} \label{1.4}$\ $ \\

Let $G$ be a graph and $o$ one particular vertex.
Consider a  random walk $(X_n)_{n\geq 0}$ on $G$ with transition
probability  $p(.,.)$ and assume there exists a reversible measure $m$ for $X$.
We use the symmetric  kernel $\mu(x,y):=m(x)p(x,y)$ to measure surfaces:
$$\forall A\subset G,\qquad \mu(\partial A)=\sum_{x\in A,y\not\in A}\mu(x,y).$$
In this setting the anchored isoperimetric inequality reads:
\begin{defi} We say $G$ satisfies the anchored isoperimetric inequality at
root $o$ with increasing function $\F$ when
for any connected $o\in A\subset G$,
\begin{equation} \label{ISW}
\frac{\mu(\partial A)}{\F(m(A))}\geq C_{IS}.
\end{equation}
``Connected'' means that one can find a discrete path in $A$
between any two points for which $p(x_i,x_{i+1})$ is positive
when $x_i,x_{i+1}$ are following points.
\end{defi}
No distance will play a role here and
the graph is not assumed to be locally finite.

We denote $\P_x$ [resp $\E_x$] the law of the walk starting from point
$x$ [resp the expectation], $\tau_A$ the exit time and $l_A$
the occupation time (which may be infinite if $X$ is not transient):
$$\tau_A=\inf\{ k\geq 0\ ;\ X_k\notin A \},\qquad
l_A={\rm card}\{ k\in \N\ ;\ X_k\in A \}.$$

\begin{theo} \label{thtl}
If $G$ satisfies (\ref{ISW}), then for any subset $A$ we have:
\begin{equation} \label{RTau}
\E_o(\tau_A) \leq 2\int_0^\infty v^A_+(s)\ {\rm d}s
\end{equation}
\begin{equation} \label{RL} and \qquad
\E_o(l_A) \leq 2\int_0^\infty v_{+A}(s)\ {\rm d}s,
\end{equation}
where $v^A,v$ are solutions of
$\begin{cases} 
v^A(0)=m(A) \\
(v^A)'= -(C_{IS}\F(v^A))^2,
\end{cases}$
and $\begin{cases} 
v(0)=+\infty \\
v'= -(C_{IS}\F(v))^2.
\end{cases}$
The truncations in indices mean
\newline $v^A_+(s)=\begin{cases}
0 \text{ if } v^A(s)\leq 0 \\
v^A(s) \text{ otherwise.}
\end{cases}$
and $v_{+A}(s)=\begin{cases}
0 \text{ if }v(s)\leq 0 \\
m(A) \text{ if }v(s)\geq m(A) \\
v(s) \text{ otherwise.}
\end{cases}$
\end{theo}
For comparison when $X$ is transient, note that
$$\int_0^\infty v^A_+(s)\ {\rm d}s
=\int_{v^{-1}(m(A))}^\infty v_+(s)\ {\rm d}s.$$
We consider usual functions $\F$ in Section \ref{ExF}. It is sometimes
useful to precise the values $\F(x)=\F(m(o))$ for $x\leq m(o)$,
which is justified in Proposition \ref{EDu}.

\section{Green functions.}

\subsection{Definitions and harmonicity.}
The submarkovian kernel of the killed random walk is $p^A(x,y)=
\begin{cases}
p(x,y) & \text{if } x \in A,   \\
0 & \text{otherwise}.
\end{cases}$
\newline Although Theorem \ref{thtl} is true for $A$ non connected, we have
in this section to assume $A$ is connected. When $X$ is transient, Green
function may be defined for the
non-killed random walk and we can consider $A=G$ (or the connected component
of $o$ if $G$ was not connected, which would have little interest). This
leads to the result for $l_A$ in next section.

The discrete Laplacian is
$$\triangle^A f=(Id -P^A)f,$$
where $P^A$ is the operator defined on functions which are zero
outside $A$ by
\begin{eqnarray*} P^Af (x)&=& \E_x(f(X_1) \ 1_{\{x,X_1\in A \} })  \\
&=& \sum_{y\in A} p^A(x,y) f(y).
\end{eqnarray*}

The Green function is
$$G^A(x,y)= \frac{1}{m(y)} \sum_{k\geq 0} \P^A_x(X_k=y).$$
In particular we denote $G^A(x)=G^A(o,x)$. Note that
$G^A(x)=0$ if $x\not\in A$.
\newline Recall that reversibility means $p(x,y)/m(y)=p(y,x)/m(x)$.
In other words $p(x,y)/m(y)$ is the precise analog of a density
kernel in $y$ starting from $x$ and is symmetric. This explains the
factor $1/m(y)$ in the definition of $G^A$ which is symmetric for
$x,y\in A$.

\begin{prop} \label{harmo}
$\triangle^A G^A = \frac{\delta_0}{m(0)}$
\end{prop}
\begin{proof} For all $x \in A $ we have :
\begin{eqnarray*} 
\triangle^A G^A (x)&=& [(Id-P^A)(G^A)] (x)\\
&=& \frac{1}{m(x)} \sum_{k\geq 0} \P^A_o(X_k=x)
-\sum_{k\geq 0} \sum_{y\in A} \frac{p^A(x,y) }{m(y)} \P^A_o(X_k=y)\\
&=& \frac{1}{m(x)} \sum_{k\geq 0} \P^A_o(X_k=x) -   \sum_{k\geq 0} \sum_{y\in A} \frac{p^A(y,x) }{m(x)} \P^A_o(X_k=y)\\
&=& \frac{1}{m(x)} \sum_{k\geq 0} \P^A_o(X_k=x) -   \sum_{k\geq 0}  \frac{1 }{m(x)} \P^A_o(X_{k+1}=x)\\
&=&\frac{ \P^A_o(X_0=x) }{m(x)}\\
&=&\frac{\delta_0(x)}{m(0)}
\end{eqnarray*}
And for $x\not\in A$, we have  $\triangle^A G^A (x) =0$.
\end{proof}

\begin{coro} \label{harmonik}$\ $ 
 $G^A$ is harmonic on $A \smallsetminus o$. As a consequence the level
sets $A_s=\{x\in A\ ; G^A(x)\geq s\}$ are connected and contain $o$.
Moreover the inward flow $through$ any $\partial A_s$ is 1 or more
generally for any $B\subset A$:
\begin{equation} \label{flo}
\sum_{(x,y)\in \partial B} \mu(x,y) \nabla_{(y,x)} G^A=1_{\{o\in B\}}.
\end{equation}
\end{coro}
The surface notations are $\partial B=\{(x,y)\ ; x\in B,y\not\in B\}$ and
$\nabla_{(y,x)}f=f(x)-f(y)$.

\begin{proof}
For all $x\in A$, Propostion \ref{harmo} may be written
$$\sum_{y\in G}  p^A(x,y) (G^A(x)-G^A(y)) =  \frac{\delta_0(x)}{m(0)}.$$
Summing over $x$ in $B$ with respect to $m$ we get
$$ \sum_{x\in B} \sum_{y\in G} m(x)  p^A(x,y) (G^A(x)-G^A(y))
=   1_{\{o\in B\}}.$$
Now the usual integration by parts becomes in this discrete summation
a cancellation of terms by symmetry when $y$ also belongs to $B$.
Only (\ref{flo}) remains.

Maximum principle and properties of level sets $A_s$ may be extracted
from this result when $o\not\in B$.
In this case the flow is $0$ so there must be an edge $x,y$ with
$G^A(y)\geq G^A(x)$. This leads to a contradiction if there was a
connected component of $A_s$ not containing $o$.
\end{proof}

\subsection{Differential inequation.}
We use a linearized version of $m(A_s)$, namely
$$u(s)=\sum_{x\in A_s,y\in G}\mu(x,y)\frac{G^A(x)-\max\{s,G^A(y)\}}
{G^A(x)-G^A(y)}.$$
For $x\in A_s$ such that $\mu(x,y)>0\Rightarrow y\in A_s$, the contribution
of $x$ is indeed $m(x)$. Furthermore $u(s)\leq m(A_s)$.
The reason for this definition is to have:
\begin{lem} Piecewise linear function $u$ has left derivative
$$u'(s)=-\sum_{(x,y)\in \partial A_s}\frac{\mu(x,y)}{\nabla_{(y,x)} G^A }.$$
\end{lem}
\begin{proof}
Variation in $s$ in the definition of $u(s)$ comes from the $y$'s
such that $G^A(y)<s$, that is $y\not\in A_s$. This is clear but note
that it uses $G^A\equiv 0$ outside $A$ and this would not be correct
for small values of $s$ and the $\tilde u$ at page \pageref{utilde} when
occupation time is considered.
\end{proof}

\begin{prop} \label{EDu} If $G$ satisfies (\ref{ISW}), then:
$$   u' \leq -(C_{\rm IS} \ \F (u) ) ^2 .  $$  
\end{prop}

\begin{proof}
Same Schwarz inequality as for Theorem \ref{thcont}:
\begin{eqnarray*}
\Big(C_{\rm IS}\F(u(s))\Big)^2 & \leq & \Big(C_{\rm IS}\F(m(A_s))\Big)^2 \\
&\leq& \mu(\partial A_s)^2 \\
&\leq& \left(\sum_{(x,y)\in A_s} \mu(x,y) \nabla_{(y,x)} G^A\right)
\left(\sum_{(x,y)\in \partial A_s}\frac{\mu(x,y)}{\nabla_{(y,x)}G^A}\right) \\
&=&-u'(s).
\end{eqnarray*}
This is of course correct when $u>0$, that is when $A_s$ is not empty and
contains $o$. It works therefore with $\F(x)=\F(m(o))$ for $x\leq m(o)$.
\end{proof}

%%%%%%%%%%%%%%%%%%%%%%%%%%%%%%%%%%%%%%%%%%%%%%%%%%%%%%%%%%%%%%%%%%%%%%%%%%%%%%%%%%%%%%%%%%%%%%%%%%%%%%%%%%%%%%%%%%%%%%%%%%%%%%%%%%%%%%%%%%%%%%%%%%%%

%%%%%%%%%%%%%%%%%%%%%%%%%%%%%%%%%%%%%%%EXIT TIME %%%%%%%%%%%%%%%%%%%%%%%%%%%%%%%%%%%%%%%%%%%%%%%%%%%%%%%%%%%%%%%%%%%%%%%%%%%%%%%%%%%%%%%%%%%%%%%

\section{Exit time}
\subsection{Upper bound}

\begin{lem} \label{lemmetime} For any set $A$ we have:
\begin{enumerate}[(i)]
\item $\E_o(\tau_A) = \sum_{x \in A} m(x) G^A(x)$,
\item $\E_o(l_A) = \sum_{x \in A} m(x) G(x)$ in the transient case.
\end{enumerate}
\end{lem}

\begin{proof} Given a path $\gamma=(\gamma_0,\gamma_1,\ldots,\gamma_n)$
from $\gamma_0=o$ to $A^c$, that is only $\gamma_n\notin A$,
we denote its probability $\P(\gamma)=p(\gamma_0,\gamma_1)\ldots
p(\gamma_{n-1},\gamma_n)$. Its length $l(\gamma)=n=\sum_{x\in A}N_x(\gamma)$
where $N_x(\gamma)$ is the number of indices $i$ such that $\gamma_i=x$.
This yields {\it (i)} since
$$\E_o(\tau_A) = \sum_\gamma l(\gamma)\P(\gamma)\quad{\rm and}\quad
G^A(x)=\frac{1}{m(x)}\sum_\gamma N_x(\gamma)\P(\gamma).$$
We adapt this argument to prove {\it (ii)}. We keep $\gamma_n\notin A$
and $\gamma_{n-1}\in A$ but we may have $\gamma_i\not\in A$ for $i<n-1$.
The probability of the path is not easy to compute but denotes
$$\P(\gamma)=\P_0(\forall i\leq n, X_i=\gamma_i
\ {\rm and}\ \forall i\geq n, X_i\not\in A).$$
We also replace the length $l(\gamma)$ by the natural occupation time
$N_A(\gamma)$.
\end{proof}

Now we could use $\sum_{x\in A}m(x)G^A(x)=\int_0^\infty m(A_s)\ {\rm d}s$.
It is a little more intricate since we have
control on $u$ which is a linearized version of $m(A_s)$.

\begin{lem} \label{lemmeintu} For any set $A$ we have:
$$\int_0^\infty u(s)\ {\rm d}s
=\sum_{x\in A,y\in G}\mu(x,y)\min\left\{G^A(x),
\frac{G^A(x)+G^A(y)}{2}\right\}.$$
\end{lem}

\begin{proof} From the definition of $u$ we just have to compute
carefully
$$\int_0^\infty\frac{G^A(x)-\max\{s,G^A(y)\}}
{G^A(x)-G^A(y)}1_{x\in A_s}\ {\rm d}s.$$
\end{proof}

We now have completed the proof of (\ref{RTau}) in
Theorem \ref{thtl}. Factor 2
in the righthand sides comes from
\begin{eqnarray*}
\sum_{x\in A} m(x)G^A(x)
&=& \sum_{x\in A,y\in G} \mu(x,y)G^A(x) \\
&\leq& 2\sum_{x\in A,y\in G}
\mu(x,y)\min\left\{G^A(x),\frac{G^A(x)+G^A(y)}{2}\right\}.
\end{eqnarray*}
As far as (\ref{RTau}) is concerned, the result first for $A$ connected
is clearly sufficient.

To prove (\ref{RL}), we first use the differential inequation
with $A=G$, that is we obtain $u(s)\leq v(s)$ for
$$u(s)=\sum_{G(x)\geq s,y\in G}\mu(x,y)\frac{G(x)-\max\{s,G(y)\}}
{G(x)-G(y)}.$$
Then we argue (here $A$ is not necessarly connected)
$$\E_o(l_A)\leq 2\sum_{x\in A,y\in G}
\mu(x,y)\min\left\{G(x),\frac{G(x)+G(y)}{2}\right\}
\leq 2\int_0^\infty \tilde u(s)\ {\rm d}s,$$
where \label{utilde}
$$\tilde u(s)=\sum_{x\in A_s,y\in G}\mu(x,y)\frac{G(x)-\max\{s,G(y)\}}
{G(x)-G(y)}.$$
It is clear that $\tilde u(s)\leq u(s)\leq v(s)$ and $\tilde u(s)
\leq m(A)$.

\subsection{Examples of $\F$ functions.}\label{ExF}
If $\F(x)=x^{1-1/d}$ as in $\mathbb{Z}^d$ then Theorem \ref{thtl}
gives
\begin{eqnarray*}
\E(\tau_A)&\leq&\frac{d}{C_{\rm IS}^2}m(A)^{2/d}\\
\text{and }
\E(l_A)&\leq&\frac{d^2}{C_{\rm IS}^2(d-2)}m(A)^{2/d} \text{ for $d>2$.}
\end{eqnarray*}
Indeed for $d>2$ the Thomassen criterium implies the transience,
see below. The computations involve
\begin{eqnarray*}
v^A(s) &=& \left(m(A)^\frac{2-d}{d}-C_{\rm IS}^2\frac{2-d}{d}s
\right)^\frac{d}{2-d} \text{ for }d\not=2, \\
v(s) &=& \left(C_{\rm IS}^2\frac{d-2}{d}s
\right)^\frac{-d}{d-2} \text{ for }d>2 \\
\text{and } v^A(s) &=& m(A)e^{-C_{\rm IS}^2s} \text{ for }d=2.
\end{eqnarray*}

If $\F(x)=x$ as in a non-amenable graph then Theorem \ref{thtl}
gives
\begin{eqnarray*}
\E(\tau_A) &\leq& \frac{1}{C_{\rm IS}^2}
\left(1+2\ln\frac{m(A)}{m(o)}\right) \\
\text{and }\E(l_A) &\leq& \frac{1}{C_{\rm IS}^2}
\left(3+2\ln\frac{m(A)}{m(o)}\right).
\end{eqnarray*}
Here we need the precision $\F(x)=m(o)$ for $x\leq m(o)$ so that 
$$\frac{1}{v^A(s)}=\frac{1}{m(A)}+C_{\rm IS}^2s$$
does not arise any issue of integration for $s\rightarrow\infty$.

We can summarize these computations in:
\begin{prop} \label{timeexit}
Let $G$ a graph satisfying a weighted anchored isoperimetric inequality with function $\F$ and anchored expansion constant $C_{IS}$ (see (\ref{ISW})). Then, there exists constants $c(d) $ and $c$ such that:
\begin{enumerate} [$\bullet$]
\item if $\F(x)=x^{1-\frac{1}{d}} \ (d\geq 3)$ we have:   $   \E_o(l_A)  \leq   c(d) \ m(A) ^ {\frac{2}{d} },$ \\
\item if $\F(x)=x^{\frac{1}{2}} \ (d=2)$ we have:   $   \E_o(\tau_A)  \leq   c(d) \ m(A) ,$ \\
\item  if $\F(x)=x$ we have:   $   \E_o(\tau_A) \leq \E_o(l_A) \leq  c \  \ln(m(A)). $ \\
\end{enumerate}
\end{prop}

\begin{rem} These inequalities are sharp. Take the particular case where $G$ satisfies a not anchored isoperimetric inequality.
\end{rem}

\begin{rem} \label{rem c(d)}Notice that the constant $c(d)$ is proportional to $1/C_{IS}^2$. There exists a constant $c_1(d)>0$ such that:
 $$c(d)=\frac{c_1(d)}{C_{IS}^2}.$$
\end{rem}
%%%%%

%%%%%%%%%%%%%%%%%%%%%%%%%%%%%%%%%%%%%%%%%%%%%%%%%%%%%%%%%%%%%%%%%%%%%%%%%%%%%%%%%%%%%%%%%%%%%%%%%%%%%%%%%%%%%%%%%%%%%%%%%%%%%%       APPLICATIONS        %%%%%%%%%%%%%%%%%%%%%%%%%%%%%%%%%%%%%%%%%%% %%%%%%%%%%%%%%%%%%%%%%%%%%%%%%%%%%%%%%%%%%%%%%%%%%%%%%%%%%%%%%%%%%%%%%%%%%%%%%%%%%%%%%%%%%%%%%%%%%
\section{Applications}
\subsection{Non degeneration for invariance principle}
As in section \ref{1.4}, assume that $X$ is  a random walk on a graph $G$ which is now supposed to be a subgraph of $\Z^d$. We suppose that $X$ admits a reversible measure $m$ satisfying $m(x) \leq b $ for all $x$ in $G$.
Under assumption of an anchored isoperimetric inequality, we will prove that if there is an invariance principle for $X$, then the diffusion constant is strictly positive. Let
 $\tilde{X}^N_k$ the renormalized random walk defined by $$\tilde{X}^N_k=\frac{1}{N} X_{kN^2}.$$
 \begin{prop}  \label{sigma}
Assume  $G$ satifies $d-$dimensionnal anchored isoperimetric inequality with constant $C_{IS}$ and that $(\tilde{X}^N_k)_k$  converges in law to a brownian motion with matrix covariance $\sigma Id$, then there exists a constant $a(d)>0$ such that  $$\sigma > \ \frac{a(d)}{b^{1/d}}  \  C_{is} .$$
In particular, $\sigma>0$. 
 \end{prop}
 
 \begin{proof}
 For a process $Z$ and a set $A$, let the correspondant exit time by $\tau_A^Z= \inf \{ k\geq 0; \ Z_k \notin A\}$. Then, by hypothesis, for all finite set $A \subset \Z^d$, we have:
 \begin{equation} \label{eq1} \lim_{N\rightarrow + \infty}  \E_0( \tau^{\tilde{X}^N }_A)  = \E_0( \tau^Y_A),\end{equation}
where $Y$ is  $d-$dimensional brownian motion with matrix covariance $\sigma^2 Id$.\\ 

Let $B(0,R)$ the ball of $\Z^d$ of radius $R$ centred at the origin. For $R>0$, we have:
\begin{enumerate} [$\bullet$]
\item
On the first hand, using the martingale $(Y_t^2 - \sigma^2t)_t$, we can prove that 
\begin{equation} \label{eq2}
\E_0( \tau^Y_{B(0,R)}) =\frac{R^2}{\sigma^2} \end{equation}
\item
On the other hand, we have:
\begin{eqnarray} \nonumber
\tau^{\tilde{X}^N }_A&=& \inf \{ k\geq 0; \ \tilde{X}^N_k  \notin A\} \\
&=& \inf \{ k\geq 0; \ X_{kN^2}  \notin N.A\} \nonumber\\
&=& \frac{1}{N^2}\inf \{ s \geq 0; \ X_{s}  \notin N.A\} \nonumber \\
&=& \frac{1}{N^2}   \tau^X_{N.A},\nonumber
\end{eqnarray}
where $N.A=\{Nx; \ x\in A\}.$\\
Hence, 
\begin{eqnarray} \nonumber 
\E_0( \tau^{\tilde{X}^N }_{B(0,R)}  )&= &\frac{1}{N^2}\E_0( \tau^{X}_{B(0,NR) } )\nonumber  \\
&\leq&  \frac{1}{N^2}  \frac{c_1(d)}{C_{IS}^2}\times  m(B(0,NR)^{2/d}   \nonumber \\
&\ &  \ \text{(by  \ proposition \ref{timeexit} and remark \ref{rem c(d)})} \nonumber\\
&\leq& \frac{c_1(d) }{C_{IS}^2}  b^{2/d} R^2.\label{eq3} \\
&\ & \text{( since $m$ is bounded by $b$.)} \nonumber
\end{eqnarray}

\end{enumerate} 

Letting $N$ goes to infinity in (\ref{eq3}), and using (\ref{eq1}) then (\ref{eq2}), we get 
$$  \frac{R^2}{\sigma^2} \leq \frac{c_1(d)}{C_{IS}^2} b^{2/d} R^2.$$
Then,
$$ \sigma^2 \geq  \frac{C_{IS}^2} {c_1(d) b^{2/d}} >0.$$
Hence, we have proved that the law of $Y$ is necessarily not degenerated
 \end{proof}
 
 \begin{rem} Let us illustrate proposition \ref{sigma} on random environments  satisfying  an  ellipticity condition. Assume that random weight verify, for all 
  edges $e$ of $\Z^d$:
$$0 \leq a \leq \omega(e) \leq b.$$
Then, for all set $A$, we have :
\begin{equation} \label{is elliptique}
\frac{\mu^{\omega} (\partial A) }{m^{\omega} (A)^{1-\frac{1}{d}}} \geq   \frac{a}{(2db)^{1-\frac{1}{d}}}  \  \frac{|\partial A|}{|A|^{1-\frac{1}{d}}}  \geq \pi_d \frac{a}{(2db)^{1-\frac{1}{d}}}
\end{equation}
where $\pi_d$ the isoperimetric constant in $\Z^d$ (see section \ref{RE} for random environment context and notations).
In particular,    (\ref{is elliptique}) is true for anchored sets. So, from  proposition \ref{sigma} we deduce if there is a invariance principle for this random walk then there exists a constant $c_d>0$ such that the diffusion constant satisfies: 
$$ \sigma \geq c_d \frac{a}{b}   .$$

 \end{rem}
%%%%%%%%%%%%%%%%%%%%%%%%%%%%%%%%%%%%%%%%%%%%%%TRANSIENCE%%%%%%%%%%%%%%%%%%%%%%%%%%%%%%%%%%%%%%%%%%%%%%%%%%%%%%%%%%%%%%%%%%%%%%%%%%%%%%%%%%%%%%%%%%%%

\subsection{Transience} We retrieve Thomassen result's  cited in the introduction. Indeed, proposition \ref{EDu} provides a new proof of the transience of the random walk under the summability assumption on $\F$ without introducing the complex construction of   dyadic  subtrees by Thomassen.
Assume \begin{equation} \label{hyp_sommabilité}
\int_1^{+\infty} \frac{1}{\F(n)^2} < +\infty,
\end{equation}
 for $\F:\R_+ \rightarrow \R_+^{\star}, \ \text{ not decreasing, with } \F(0)=0$ and let us prove transience with the  help of Proposition \ref{EDu}. 

Let $A$ a connected subset of $G$ containing the origin  and consider random walk killed whenever it leaves $A$ and the associated Green function $G_A$.
Integrating the differential equation of Proposition \ref{EDu} between time 0 and t gives:
\begin{equation}
\label{edu intégrée}  
\int_{u(t)}^{u(0)} \frac{ds}{\F(s)^2}  \geq C_{IS}^2 \ t.
\end{equation}
$ \int_{1}^{u(0)} \frac{ds}{\F(s)^2}  $ is bounded by a constant independant of $A$. Indeed, thanks to  hypothesis (\ref{hyp_sommabilité}), for all subset $A$ we have:
 $ \int_{1}^{u(0)} \frac{ds}{\F(s)^2}  = \int_{1}^{m(A)} \frac{1}{\F(s)^2} ds \leq 
\int_{1}^{+\infty} \frac{ds}{\F(s)^2}  < +\infty.$
So for large enough t which depends only on $ C_{IS} $ and $\F$, inequality (\ref{edu intégrée}) turns into: 
$$ \int_{u(t)}^{1} \frac{ds}{\F(s)^2}  \geq  \frac{1}{2} C_{IS}^2 \ t.$$
Then, we deduce that: $$\lim_{t\rightarrow +\infty} u(t)=0 \ \text{ uniformly in } A. $$ In particular,
there exists $t_0$ independant of $A$ such that for all $t\geq t_0, \ u(t) < \inf_G m$.
Therefore by definition of $u$ we get that for all set $A$, $G_A\leq t_0$. Now we can make $A$ growing and finally we deduce that $G < +\infty$ so the walk is transient.
%%%%%%%%%%%%%%%%%%%%%%%%%%%%%%%%%%%%%%%%%%%%%%%%%%%%%%%% SPEED %%%%%%%%%%%%%%%%%%%%%%%%%%%%%%%%%%%%%%%%%%%%%%%%%%%%%%%%%%%%%%%%%%%%%%%%%%%%%%%%%%%%%

\subsection{Speed}
When $\F=id$, the upper bound of the exit time gives us that the speed of the random walk is  positive. We retrieve a weak version of Virag's result. We assume in this subsection that the graph has uniformly localy bounded valency.
Let $d(a,b)$ denote the graph distance between point $a$ and $b$.
\begin{prop}
Let $G$ be a graph satisfying (\ref{IS}) with $\F=id$ and let
$(X_n)_n$ be a simple random walk on $G$. Then we have:
$$  \P\left(  \   \lim_{n}   \frac{d(o,X_n) }{n} =0\right) =0. $$
\end{prop}
\begin{proof}
Assume there exists $\epsilon >0 $ such that  $\P( \lim_n \frac{d(o,X_n)}{n} =0) >\epsilon $.  So, we have:
$$\forall \alpha >0 \ \ \ \P(\exists N_{\alpha} \  \forall n\geq N_{\alpha}  \  \ \  \frac{d(o,X_n)}{n} \leq \alpha) >\epsilon  $$
By considering the event $E_q=\{ \exists N_{\alpha} <q, \ \forall n\geq N_{\alpha}  \  \ \  \frac{d(o,X_n)}{n} \leq \alpha \}$ and by continuity of measure $\P$, we get: 
\begin{eqnarray}\label{Nalpha}  \exists N_{\alpha}\geq 0 \ \ \ \P\left(\forall n\geq N_{\alpha}  \ \ \ \frac{d(o,X_n)}{n} \leq \alpha\right) > \frac{\epsilon}{2}.  \end{eqnarray}

Take now $R>0$, we have:

$$\P\left( \forall n \in [N_{\alpha }; \frac{R}{\alpha}] \ \ \ d(o,X_n) < \alpha n     \right) > \frac{\epsilon}{2 }. $$
On this event we have:   $l_{B(o,R) } \geq \frac{R}{\alpha} - N_{\alpha}$, where $l_A$ is the local time of $X$ in the set $A$, which is well defined in this case since when $\F=id$ the walk is transient by Thomassen result.
Therefore, by using (\ref{Nalpha}), we get \begin{eqnarray} \label{einf} \E_o(l_{B(o,R) } ) \geq  \frac{\epsilon}{2} \left(\frac{R}{\alpha} -N_{\alpha}\right). \end{eqnarray}
By Proposition \ref{timeexit} and since strong anchored isoperimetric inequality implies a subexponential volume growth, there exists $c>0$ such that:
 \begin{eqnarray} \label{esup}
\E_o(l_{B(o,R) } ) \leq \ln ( |B(o,R)|) \leq  c R  
 \end{eqnarray}

Choose now $\alpha$ such that $ \frac{\epsilon}{2\alpha} >c$. Gathering (\ref{einf})and (\ref{esup}), we get:
$$\frac{\epsilon}{2} \left(\frac{R}{\alpha} -N_{\alpha}\right) \leq c R$$
  Letting $R$ goes to infinity in this last expression, we get a contradiction.
\end{proof}

%%%%%%%%%%%%%%%%%%%%%%%%%%%%%%%%%%%%%%%%%%%%%%%%%%%%%%%%%%%%%%%%%%%%%%%%%%%%%%%%%%%%%%%%%
%%%%%%%%%%%%%%%%%%%%%%%%%       RANDOM ENVIRONMENTS      %%%%%%%%%%%%%%%%%%%%%%%%%%%%%%%%% %%%%%%%%%%%%%%%%%%%%%%%%%%%%%%%%%%%%%%%%%%%%%%%%%%%%%%%%%%%%%%%%%%%%%%%%%%%%%%%%%%%%%%%%%

\section{Random environment on $\Z^d$, including supercritical percolation}  \label{RE}
We consider discrete time, nearest-neighbor random walks
among random (i.i.d.) conductances on $\Z^d, d\geq 2$.
Our model will include super-critical percolation since
conductances may be null, we do not require that
conductances be bounded, just that they are exponentially
integrable.

After a presentation of random environment in the first
subsection, we prove an isoperimetric inequality for big sets
in the second part, which leads to occupation time estimate
for big sets in the third subsection.

\subsection{Setting: super-critical exponentially integrable
random environment} 
Consider graph $\L^d=(\Z^d,E_d)$ where $E_d$ contains
non-oriented nearest-neighbor pairs. We write $x\sim y $
if $\{x,y\} \in E_d$.  
An environment is a random function $\mu^\omega : E_d \rightarrow
[0;+\infty[$. It is implicit in the definition of $E_d$ that it is
symetric. The value $\mu^\omega(x,y)$ is called the conductance
of edge ${x,y}$. To lighten the notations, we will sometimes write $\mu$ instead of $\mu^{\omega}$ when there is no ambiguity.

Let $\mathbb{Q}$ be a product probability measure on
$[0;+\infty[^{E_d}$.
A walker or an electric current can cross only edges with
strictly positive conductances. So we call cluster a connected
component of the graph $(\Z^d, \{e\in E_d\ ; \ \omega(e)>0 \})$
and we use $\mathbb{Q}$-connectedness refering to this graph.
In fact $Q$ induces a Bernoulli percolation $P_Q$ of parameter
$q=\mathbb{Q}(\mu(e)>0)$ (here and in the following, $e$ is any
edge since $\mathbb{Q}$ is a product measure).
We assume $q>p_c$ critical parameter of edge percolation on $\Z^d$
\begin{defi} For $q>p_c$, the law $\mathbb{Q}$ is said to be a  super-critical exponentially
integrable random environment if  there exists $\beta>0$ such that
$$\E_{\mathbb{Q}}(\exp(\beta\mu(e)))<\infty.$$
\end{defi}

For quenched result on the random walk, we consider measure $\P_0=
\P_{0,\mu^\omega}(\cdot\mid \C_0\text{ infinite})$. That is, we start
the random walk from the origin $0$ of $\Z^d$, $\mu$ induces $m$ and
$p(\cdot,\cdot)$ so that we are in the setting defined in Section
\ref{1.4}, and we assume the cluster $\C_0$ of $0$ is infinite.

 %%%%%%%%%%%%%%%%%%%%%%%%%%%%%%%%%%%%%%%%%%%%%%%%%%%%%%%%%%%%%%%%%%%%%%%%%%%%%%%%%%%%%%%%
%%%%%%%%%%%%%           IS ancrée sur environement        %%%%%%%%%%%%%%%%%%%%%%%%%%%%%%% %%%%%%%%%%%%%%%%%%%%%%%%%%%%%%%%%%%%%%%%%%%%%%%%%%%%%%%%%%%%%%%%%%%%%%%%%%%%%%%%%%%%%%%%%

\subsection{ (Anchored) isoperimetric inequality}

We need an anchored isoperimetric inequality with respect to random
weight $\mu^\omega$. Differents forms of strong isoperimetric inequality
have been established by many authors (see \cite{pierre}, \cite{KL}, \cite{pete} and \cite{bisk})
in the percolation context.

We may only have a control for big sets.  
The form which seems adapted to our exit time results
is the following:
\begin{prop} \label{is ancrée RE}  Let  $\mathbb{Q}$ a super-critical
exponentially integrable random environment on $\Z^d$.
\newline There exist $\beta_0(\mathbb{Q})>0$ and a random integer $N_0(\omega)$
such that,
\newline for all $\mathbb{Q}$-connected sets $A\subset\Z^d$ containing 0,
 \begin{eqnarray} \label{bbetac}
 |A|\geq N_0(\omega) \qquad
 \Longrightarrow \qquad
 \frac{\mu^\omega(\partial A)}{m^\omega(A)^{1-1/d}} \geq \beta_0.
 \end{eqnarray}
\end{prop}

%%%%%%%%%%%%%%%%%%%%%%%%%%%%%%%%%%%%%%%%%%%%%%%%%%%%%%%%%%%%%%%%%%%%%%%%%%%%%%%%%%%%%%%
%%%%%%%%%%%%%%%%%%%%%%%%%%%%%%     LEMME NEW     %%%%%%%%%%%%%%%%%%%%%%%%%%%%ù %%%%%%%%%%%%%%%%%%%%%%%%%%%%%%%%%%%%%%%%%%%%%%%%%%%%%%%%%%%%%%%%%%%%%%%%%%%%%%%%%%%%%%%%%%%%%%%%%%%%

In order to use non-weighted isoperimetric inequality
$|\partial A|/|A|^{1-1/d}\geq C_d$,
the first point in the proof is to control $m(A)$ with $|A|$.
The difficulty is it has to be done for any $A$. But a standard
exponential Bienaymee Tchebytchef inequality works.

\begin{lem}\label{lemmeEBT1} There exist $\beta_1(\mathbb{Q})>0$
and a random integer $N_1$ such that,
\newline for all $\mathbb{Q}$-connected sets $A\subset\Z^d$ containing 0,
$$|A|\geq N_1 \qquad\Longrightarrow \qquad
|A|\geq \beta_1 m(A).$$
\end{lem}

\begin{proof}
Denote $A_n=\{A\subset\Z^d\ ;\ A\text{ is }\mathbb{Q}
\text{-connected };\ 0\in A
\text{ and } |A|=n\}$ and $\A_A=\{|A|<\beta_1 m(A)\}$.
\newline We should prove that $\displaystyle \A_n=
\bigcup_{A\in A_n}\A_A$ may not occur infinitely often.

By exponential Bienaymee Tchebytchef inequality,
$$\mathbb{Q}(\A_A)\leq e^{-\lambda |A|}\E e^{\lambda\beta_1 m(A)}
\leq e^{-\lambda n} (\E e^{2\lambda\beta_1\mu(e)})^{dn},$$
if $A\in A_n$, since then $m(A)$ is a summation of $dn$
variables $\mu$ with a factor $2$ or $1$.

Since $A$ is connected and contains $0$,
there exists $\alpha(d)$ such that
$|A_n|\leq e^{\alpha n}$. This may be shown
by constructing possible $A$'s starting from $0$.
This leads to
$$\mathbb{Q}(\A_n)\leq e^{\alpha n}
e^{-\lambda n} (\E e^{2\lambda\beta_1\mu(e)})^{dn}.$$
To use Borel-Cantelli lemma and finish the proof, we
should find coefficients $\lambda$ and $\beta_1$ such that
the right-hand side is summable in $n$.
\newline First step is $\lambda=2\alpha$.
\newline Then $\beta_1$ is chosen such that
$(\E e^{2\lambda\beta_1\mu(e)})^d<e^\alpha$, which is possible since
exponential integrability implies that $\E e^{\beta\mu(e)}$ tends
to $1$ when $\beta$ tends to $0$.
\end{proof}

The case of $\partial A$ is more involved. This frontier may
be reduced by small values of $\mu$ or by null values, that is
percolation. The exponential Bienaymee Tchebytchef argument
works when percolation parameter $q$ is close to $1$, this is
Lemma \ref{lemmeEBT2}.

After this lemma, we shall use renormalization. So $F_2$ in Lemma
\ref{lemmeEBT2}
will be a subset of boxes around $\partial A$. We denote
$3^d$-connectivity to include boxes sharing only one corner,
precisely $x\sim y$ if $|x-y|_\infty\leq 1$. We consider Bernoulli
site percolation $\P_{p_2}$ (bad or good box) of parameter $p_2$.
The subset $\tilde F_2\subset F_2$ designs open sites of $F_2$.

A small additional difficulty comparing to Lemma \ref{lemmeEBT1}
is the loss of anchor $0$, since $F_2$ is rather $\partial A$.
We need that $F_2$ stays in some big box and we will use the scale
$(\ln n)^{3/2}$ for $F_2$ and $n$ for the big box. 
% like in piatnitski and Mathieu PP
We note ${\rm Box}(n)=\Z^d\cap[-n,n]^d$.

\begin{lem}\label{lemmeEBT2}
There exist $\beta_2>0$, $p_2<1$ and a random integer $N_2$ such that,
\newline  $\P_{p_2}$ almost surely for all $n\geq N_2$ and
$F_2\subset{\rm Box}(n)$ $3^d$-connected,
$$|F_2|\geq(\ln n)^{3/2}\qquad\Longrightarrow\qquad
|\tilde F_2|\geq \beta_2|F_2|.$$
\end{lem}

\begin{proof}
Denote $A_n=\{F_2\subset\{-n,\ldots,+n\}^d\ ;
\ F_2\text{ is }3^d\text{-connected }
\text{ and } |F_2|\geq(\ln n)^{3/2}\}$,
$A_{n,m}=\{F_2\in A_n\ ;
\ |F_2|=m\}$
and $\A_{F_2}=\{|\tilde F_2|<\beta_2 |F_2|\}$.
\newline We should prove that $\displaystyle \A_n=
\bigcup_{F_2\in A_n}\A_{F_2}$ may not occur infinitely often.

By exponential Bienaymee Tchebytchef inequality,
$$\P_{p_2}(\A_{F_2})\leq e^{\lambda\beta_2|F_2|}\E e^{-\lambda|\tilde F_2|}
= e^{\lambda\beta_2 m}(1-p_2+p_2e^{-\lambda})^m,$$
if $F_2\in A_{n,m}$.

Since $F_2\subset\{-n,\ldots,+n\}^d$ is connected,
starting from any point of $F_2$, we may show
$$|A_{n,m}|\leq(2n+1)^d e^{\alpha'm}.$$
Since $m\geq(\ln n)^{3/2}$, this may be simplified
$|A_{n,m}|\leq e^{\alpha m}$ for some $\alpha$.
This leads to
$$\P_{p_2}(\A_n)\leq \sum_{m\geq(\ln n)^{3/2}} e^{\alpha m}
e^{\lambda\beta_2 m}(1-p_2+p_2e^{-\lambda})^m.$$
To use Borel-Cantelli lemma and finish the proof, we
first choose $\lambda>\alpha$.
\newline Then $p_2$ is chosen such that
$\displaystyle 1>e^\alpha(1-p_2+p_2e^{-\lambda})\underset
{p_2\rightarrow 1}\longrightarrow e^\alpha e^{-\lambda}$.
\newline Finally $\beta_2$ is chosen such that the
right-hand side is bounded, for some $\gamma>0$, by
$$\sum_{m\geq(\ln n)^{3/2}} e^{-\gamma m}
\leq\frac{1}{1-e^{-\gamma}}e^{-\gamma(\ln n)^{3/2}},$$
which is summable in $n$.
\end{proof}

Since in Proposition \ref{is ancrée RE} we only assume
$q>p_c$, we will use renormalization, namely
Proposition 2.1 from \cite{antal}.

For length $L$, $\Z^d$ is parcelled into boxes
$B_i(L)=\tau_{i(2L+1)}{\rm Box}(L)$ where $\tau$ designs
translation. We define a percolation $\P_L$ on the $i$'s
induced by Bernoulli percolation $\P_q$. The site $i$ is
open if the box $B_i$ is ``good'' in the following sense~:
there is a unique crossing cluster $\C\subset B'_i
=\tau_{i(2L+1)}{\rm Box}(5L/4)$, which means
that for all $d$ directions, $\C$ joins the two faces of $B_i$.
Furthermore any open path in $B'_i$ of length larger than
$L/10$ is connected to $\C$ in $B'_i$, and $\C$ is crossing for all subbox $B\subset B'_i$ of size length larger than $L/10$. A lot more is demanded
in \cite{antal}, but this is sufficient here.

\begin{prop} \label{renorm}
For all $q>p_c$ there exists $L(q)$ such that
$$\P_L\geq\P_{p_2}.$$
\end{prop}

This stochastic inequality should be understood for increasing
events like $\A_{F_2}$ of Lemme \ref{lemmeEBT2}. In fact we will
apply this lemma to $\P_L$.

We finish with the exponential Bienaymee Tchebytchef technique
for small values of $\mu$. The connectivity argument gets still
more technical since we will consider some $\tilde F_3$ already
reduced by percolation, and this $F_3$ will contain one edge for
each big box. Thus we define $L$-connectivity of edges by the existence
of a path of length $10L$ between them. We denote $E^d(n)$ the set
of edges between sites of ${\rm Box}(n)$ and consider $\mathbb{Q}$,
denoting $\tilde F_3=\{e\in F_3\ ;\ \mu(e)>0\}$.

\begin{lem}\label{lemmeEBT3}
Let $L=L(p_2)$ and $\beta_2$ given by Lemma \ref{lemmeEBT2}.
There exist $\beta_3>0$ and a random integer $N_3$ such that,
\newline for all $n\geq N_3$ and $F_3\subset E^d(n)$ $L$-connected,
$$|F_3|\geq(\ln n)^{3/2}\text{ and }|\tilde F_3|\geq\beta_2 |F_3|
\qquad\Longrightarrow\qquad
\mu(F_3)\geq \beta_3|F_3|.$$
\end{lem}

\begin{proof}
Denote $A_n=\{(F_3,\tilde F_3)\ ;
\ \tilde F_3\subset F_3\subset E^d(n),
\ F_3\text{ is }L\text{-connected}
\text{ and } |F_3|\geq(\ln n)^{3/2}\}$,
$A_{n,m}=\{(F_3,\tilde F_3)\in A_n\ ;\ |F_3|=m\}$
and $\A_{F_3,\tilde F_3}=\{\mu_{|\tilde F_3}>0 \text{ and }
\mu(F_3)<\beta_3 |F_3|\}$.
\newline As usual we should prove that $\displaystyle \A_n=
\bigcup_{(F_3,\tilde F_3)\in A_n}\A_{F_3,\tilde F_3}$ may not occur
infinitely often.

By exponential Bienaymee Tchebytchef inequality,
$$\mathbb{Q}(\A_{F_3,\tilde F_3})\leq e^{\lambda\beta_3|F_3|}
\E (e^{-\lambda\mu(F_3)}\mid \mu_{|\tilde F_3}>0) \leq e^{\lambda\beta_3 m}
(\E e^{-\lambda\mu(e)}\mid \mu(e)>0)^{\beta_2 m},$$
if $(F_3,\tilde F_3)\in A_{n,m}$.

For some $\alpha$, $|A_{n,m}|\leq e^{\alpha m}$. Indeed the
choice of subset $\tilde F_3$ may be bounded by $2^m$.
This leads to
$$\mathbb{Q}(\A_n)\leq \sum_{m\geq(\ln n)^{3/2}}
e^{\alpha m} e^{\lambda\beta_3 m}
(\E e^{-\lambda\mu(e)}\mid \mu(e)>0)^{\beta_2 m}.$$
Now we first choose $\lambda$ big enough so that
$(\E e^{-\lambda\mu(e)}\mid \mu(e)>0)^{\beta_2}<e^{-\alpha}$,
and then $\beta_3$ small enough to estimate the
right-hand side by some $\sum e^{-\gamma m}$.
\end{proof}

{\bf Proof of Proposition \ref{is ancrée RE}.} \\
Let $L=L(p_2)$ given by Lemma \ref{lemmeEBT2} and Proposition
\ref{renorm}. We cover the outer boundary of $A$ by boxes $B_i(L)$
the following way: $i\in F_2$ means that $A\cap B_i\not=\emptyset$
and at least one face of $B_i$ stays in the infinite component of
$\Z^d\setminus A$.

Under $\P_L$, for any ``good'' box $B_i$ where $i\in F_2$,
we show how one can find
an edge $e_i\in\partial A\cap B'_i$ such that $\mu(e)>0$, (see figure 1 below).We start
from $x\in A\cap B_i$, its component in $A\cap B_i$ is big enough
so that $x$ belongs to the crossing cluster (either take a path
from x in $A$ to  outside $B'_i$, either take $A$ big enough in the very
special case when $A\subset B'_i$, use $N_0$ below).
Now take a path in the crossing cluster
from $x$ to the face justifying $i\in F_2$. On this path, $\mu>0$, it
begins in $A$ and leaves it at some $e_i$.

We collect these $e_i$'s in $F_3$. We want that $|F_3|=|F_2|$ and that
$F_3$ is $L$-connected, so we may add any edge of $\partial A
\cap B_i$ for ``bad'' boxes. Also we may have to take any edge in
the special case when two close ``good'' boxes give the same $e_i$.
Anyway we can have $|\tilde F_3|\geq|\tilde F_2|/2^d$.
If we can (by fixing $N_0$ below) apply
the three lemmas and Proposition \ref{renorm}, we will have
$$\mu(\partial A)\geq\beta_3|F_3|\geq\frac{\beta_3}{2d(2L+1)2^d}
C_d|A|^{1-1/d}\geq\frac{\beta_3\beta_1^{1-1/d}C_d}{2d(2L+1)2^d}
m(A)^{1-1/d}.$$
We have used successively~: Lemma \ref{lemmeEBT3} since
$F_3\subset\partial A$, non-weighted isometric inequality for $A$
which may be surrounded by all faces of the $B_i$'s ($i\in F_2$)
and Lemma \ref{lemmeEBT1}.

There remains to check that a good choice of $N_0$ function of $N_1$,
$N_2$ and $N_3$ ensures that $|A|\geq N_0$ implies all technical
conditions:
\begin{itemize}
\item $|A|\geq N_1$ for Lemma \ref{lemmeEBT1} is easy.
\item $n\geq N_2$, $F_2\subset {\rm Box}(n)$
and $|F_2|\geq(\ln n)^{3/2}$ for Lemma \ref{lemmeEBT2}.
\item $n\geq N_3$ and $F_3\subset E^d(n)$ for Lemma \ref{lemmeEBT3}.
\end{itemize}
For the last point, we need $A\subset {\rm Box}(n-9L/4)$
since the $e_i$'s may be chosen in $B'_i$. Take $n$ the smallest
integer such that this is true. We have $f(n)\leq|A|\leq g(n)$ for
some functions $f$ and $g$ which behave asymptotically like $n$
and $n^d$.

Thus a good choice of $N_0$ makes sure that $g(n)\geq|A|\geq N_0$
implies $n\geq N_2$ and $n\geq N_3$.
Since the union of boxes of $F_2$ surround $A$, we have
$$(2L+1)^d|F_2|\geq C_d f(n)^{1-1/d}.$$
It will ensure $|F_2|\geq(\ln n)^{3/2}$ if $N_0$ is big enough.

%\begin{figure}[center]
%   \caption{\label{fig1}  }
%\includegraphics[scale=0.35]{isdessin.pdf}
%\end{figure}
\begin{rem} In \cite{bouk}, it is proved that we can build  environments where the return probability is greater than $1/n^2$. By our proposition \ref{is ancrée RE}, the d-dimensional anchored isoperimetric inequality is satisfied on these environments and so in dimension higher than $4$, no one can hope to prove that in this case, the return probability is in $1/n^{d/2}.$ 
\end{rem}

\begin{rem} \label{A petit} Let $\omega$ a fixed environment and $N_0(\omega)$ as in Proposition \ref{is ancrée RE}. Since $Q(\omega>0)=1$, there is a finite number of sets $B$ containing $0$ and satisfying $m^{\omega} (B) \leq N_0(\omega)$. Thus for a set $A$ such that $m^{\omega}(A)\leq N_0(\omega)$, we can have 
$$\mu^{\omega}(\partial A ) \geq c_{\omega}:= \min \{  \sum_{e\in \partial B} \omega(e); \ 0 \ni B \ such \ that \ m^{\omega} (B) \leq N_0(\omega) \}>0. $$ This  can be re written as well as follow:
\begin{equation}\label{is A petit}
\mu^{\omega}(\partial A )=
\sum_{e\in \partial A}
 \omega(e) \geq  
 \beta_{\omega} \  
 m^{\omega}( A)^{1-1/d},
\end{equation} with
$\beta_{\omega}=c_{\omega}/N_0(\omega)^{1-1/d}, $ constant which depends on $\omega$.  
\end{rem}

%%%%%%%%%%%%%%%%%%%%%%%%%%%%%%%%%%%%%%%%%%%%%%%%%%%%%%%%%%%%%%%%%%%%%%%%%%%%%%%%%%%%
%%%%%%%%  Temps d'occupation sur environments aléatoire  %%%%%%%%%%%%%%%%%%%%%%%%%%% %%%%%%%%%%%%%%%%%%%%%%%%%%%%%%%%%%%%%%%%%%%%%%%%%%%%%%%%%%%%%%%%%%%%%%%%%%%%%%%%%%%%

\subsection{Upper bound for the occupation time ($d\geq 3)$ and for the exit time ($d=2$)}
We apply result of Theorem \ref{thtl}  in the particular case of random walk in random environment with $\mathbb{Q}$ a super-critical exponentially
integrable random environment. We get,

\begin{prop} \label{propRE}   There exists constant  $C= C(Q,d)$ such that  $Q$ a.s. for all environment $\omega$,  $Q$ super-critical exponentially integrable environment:\\
  for any  connected  subset $B$  which contains the origin and with volume $|B|$ large enough, 

\begin{numcases}{}\label{tempsperco1}
 \E_0(l_B)  \leq C m^{\omega}(B)   \ \ \ \text{ in dimension } d\geq 3 \\
             \label{tempsperco2}  
 \E_0(\tau_B)  \leq C m^{\omega}(B) \ \ \ \text{ for dimension } d=2
\end{numcases}
\end{prop}

$\ $ 
\\
Let  $B\subset \Z^d$ connected and which contains the origin.  We are going to estimate $\E(\tau_B)$ (or $\E(l_B)$ in transient case).

%%%%%%%%%%%%%%%%%%%%%%%%%%%%%%%%%%%%%%%%%%%%%%%%%%%%%%%%%%%%%%%%%%%%%%%%%%%
%%%%%%%%%%            d larger than 3  %%%%%%%%%%%%%%%%%%%%%%%%%%%%%%%%%%%% %%%%%%%%%%%%%%%%%%%%%%%%%%%%%%%%%%%%%%%%%%%%%%%%%%%%%%%%%%%%%%%%%%%%%%%%%%%

\begin{proof} $\ $ \\
(i) {\bf{case $d\geq 3$}}

By our isoperimetic inequality (Proposition \ref{is ancrée RE} 
%and remark \ref{A petit}
),  and by result of Thomassen, we deduce that the walk is transient. So we can  deal with  $G$ the  whole Green fonction. 
For $t\geq 0$ we  let $$u(t)= m^{\omega} (\{  x\in B; \ G(0,x) \geq t \}).$$ By Proposition \ref{EDu} and thanks to inequality (\ref{bbetac}),  function $u$ satisfies:
$$   \begin{cases}
u(0)=m^{\omega}(B)   \\
 u'\leq  - (\beta_0 \ u^{1-\frac{1}{d}})^2,\; \text{ until $\#\{ x\in B;\ G(0,x)\geq t \}\geq N_0(\omega)$}\\
\end{cases}   $$

Assume that $|B|\geq N_0(\omega)$. Solving this differential equation, we get:
\begin{eqnarray} \label{u perco}
 u(t) \leq 
\left[ \ \frac{d -2}{d} \beta_0^2t + m^{\omega}(B)^{\frac{2}{d}-1  }  \ \right]^{\frac{d}{2-d}} & \text{if } t\leq t_0  
\end{eqnarray}
with $t_0$ such that $$ \#\{ x\in B;\ G(0,x)\geq t_0 \}\geq N_0(\omega).$$
 Now, Corollary 3.2 gives us the expectation of the occupation time. We have:
$$ \E_0(l_B) = \int_0^{+\infty } u(s) \ ds$$
We  split into two parts the computation of  this integral.
First, we have:
 $$\int_{0}^{t_0} u(s) \ ds \leq  \frac{d}{ 2\beta_0^2} \left[ m^{\omega}(B)^{\frac{2}{d}} - M_0(\omega)^{\frac{ 2}{d}  }  \right],$$
 with $M_0(\omega):=m^{\omega}(\{ x\in B;\ G(0,x)\geq t_0 \})$.
 Secondly we have to deal with the term $\int_{t_0}^{+\infty} u(s) \ ds $
 \begin{eqnarray*}
 \int_{t_0}^{+\infty}u(t) \ dt &=& \int_{t_0}^{+\infty} m( \{ x\in B; \  G(x)\geq t \}) \ dt \\
 &=& \int_{t_0}^{G(0)} m( \{ x\in B; \  G(x)\geq t \}) \ dt \\
 &\leq&  (G(0)-t_0)m( \{ x\in B; \  G(x)\geq t_0 \})\\
 &=&  (G(0)-t_0)M_0(\omega)\\
 &\leq&  G(0) M_0(\omega)
\end{eqnarray*}  
Gathering the two previous computations, we  get:
$$ \E_0(l_B) \leq \frac{d}{ 2\beta_0^2} \left[ m^{\omega}(B)^{\frac{2}{d}} - M_0(\omega)^{\frac{ 2}{d}  }  \right] + G(0) M_0(\omega).$$
Finally, we have proved that there exists $C>0$ such that for $Q$ a.s. environment $\omega$, there exists $N_{\omega}\in \N$ such that for any  connected  subset $B$ which contains the origin with, $m^{\omega}(B)\geq   N_{\omega}$ then $\E_0(l_B)  \leq C m^{\omega}(B)^{2/d}$.

%%%%%%%%%%%%%%%%%%%%%%%%%%%%%%%%%%%%%%%%%%%%%%%%%%%%%%%%%%%%%%%%%
%%%%%%%%%%%%%%%%%        d = 2    %%%%%%%%%%%%%%%%%%%%%%%%%%%%%%%
%%%%%%%%%%%%%%%%%%%%%%%%%%%%%%%%%%%%%%%%%%%%%%%%%%%%%%%%%%%%%%%%%%

(ii) {\bf{case $d=2$}}

The same kind of  arguments gives the bound in the dimension two replacing the occupation time by the exit time. In recurrence case, for  $t\geq 0$ we  let $$u(t)= m^{\omega} (\{  x\in B; \ G^B(0,x) \geq t \}),$$ where $G^B$ is the Green function of the random walk killed outside $B$.
Once again, we use  isoperimetric inequality of Proposition \ref{is ancrée RE} (equation (\ref{bbetac})).
% and moreover we use remark \ref{A petit} (equation (\ref{is A petit})).
By Proposition \ref{EDu}, for set $B$ such that $|B|\geq N_0(\omega)$, we get:
\begin{eqnarray} \label{u d=2}
 u(t) \leq  
 %\begin{cases}
 m^{\omega}(B)e^{-\beta_0^2t} & \text{until } t\leq t_0,
 % N_0(\omega)e^{-\beta_{\omega}^2(t-t_0)}& \text{if } t>t_0.
%\end{cases}
\end{eqnarray}
with $t_0$ such that $\#\{x\in B; G^B(0,x)\geq t_0\} \geq N_0(\omega). $
If we let $M_0(\omega)=m^{\omega}(\{x\in B; G^B(0,x)\geq t_0\}),$ then $t_0$ satisfies $t_0=\frac{1}{ \beta_0^2 } \ln \left( \frac{m^{\omega}(B)}{ M_0(\omega)}\right)    $
Then, 
 \begin{eqnarray*}
 \E_0(\tau_B) &=& \int_{0}^{+\infty}u(s) \ ds \\
 &\leq& m^{\omega}(B) 
 [ \frac{1}{\beta_0^2} 
 +\frac{N_0(\omega)}{   m^{\omega}(B)}
% (\frac{1}{\beta_{\omega}^2 } 
 %-\frac{1}{\beta_0^2})
 C(\omega)] 
\end{eqnarray*} 
%Finaly, we get the same proposition.
\end{proof}

%%%%%%%%%%%%%%%%%%%%%%%%%%%%%%%%%%%%%%%%%%%%%%%%%%%%%%%%%%%%%%%%%%%%%%%%%%%%%%%
%%%%%%%%%%%%%%%%%%%%    REMARQUE  PERCOLATION    %%%%%%%%%%%%%%%%%%%%%%%%%%%%%% %%%%%%%%%%%%%%%%%%%%%%%%M%%%%%%%%%%%%%%%%%%%%%%%%%%%%%%%%%%%%%%%%%%%%%%%%%%%%%% 
\subsection{Percolation case}
Percolation is a particular case of $Q$ super-critical exponentially integrable random environment. So, theses results hold for percolation super-critical cluster.
\begin{prop} \label{tau sur perco}  Let $p> p_c(d)$ and $d\geq 2$. There exists constant $C= C(p,d)$ such that  $Q$ a.s. on the event  $\{\#\C=+\infty\}$:\\
  for any  connected  subset $B$ of $\C$ which contains the origin and with   volume large enough, 
\begin{equation} \label{tempsperco}
\left\{
          \begin{array}{ll}
 \E_0(l_B)  \leq C |B|^{2/d} \ \ \  \text{if } d\geq 3,\\
 \E_0(\tau_B)  \leq C |B|^{2/d}  \ \ \     \text{if } d \geq 2.
\end{array}
        \right.
\end{equation} 
\end{prop}

\begin{rem}  These estimates have the right behaviour, since we retrieve a consequence of results of Barlow. Indeed, in \cite{barlow} it is proved that:
\begin{theo}There exists $\Omega_1$ with $Q(\Omega_1)=1$ and  random variables $ Sx; x\in\Z^d$  such that for each $x\in\C$ and for all $\omega\in\Omega_1, \  S_x(\omega) < \infty $ and there exist constants $ci = c_i(d; p)>0$ such that
for all $x, y \in \C$ and $t\geq 1$ with
$$ k\geq S_x(\omega) \vee |x-y|_1 $$ 
the transition density  $\P_x(X_k=y)$ of $X$ satisfies:
$$\nu(y) c_1 k^{-d/2} e^{-c_2 \frac{|x-y|_1^2}{k}} \leq \P_x(X_k=y)\leq \nu(y) c_3 k^{-d/2} e^{\frac{-c_4 |x-y|_1^2}{k}}.$$
\end{theo}
\end{rem}
Let $x_0 \in \C$ and let $B\subset \C$ connected which contains the point $x_0$.  First, for all $k\geq 0$ we can write:
$$\P_{x_0}(\tau_B>k) \leq \sum_{y\in B} \P_{x_0}(X_k=y).$$
With the help of the previous Theorem (we keep the same notation),  there exists a constant $c>0 $ such that $Q$ a.s.  for all $ y\in B $ and for all $k\geq S_{x_0}(\omega) \vee |x_0-y|_1$, we have :
$$ \P_{x_0}(X_k=y)\leq c \nu(y)  k^{-d/2}.$$
Hence,$$ \P_{x_0}(\tau_B>k) \leq  c \nu(B)k^{-d/2}. $$
Let  $k_0= (2c . \nu(B))^{2/d}$ and fix  the environement $\omega$. For $B$ large enough, the condition $k_0\geq S_{x_0}(\omega) \vee |x_0-y|_1$  is satisfied (once again the size from which the condition is satisfied depends on the point $x_0$ and $\omega$).

Then, for all $x_0$ and for $B$ large enough, connected and which contains $x_0$, we have :
$$ \P_{x_0}(\tau_B>k_0) \leq 1/2. $$
So, for all $i\geq 0,$ $$ \P_{x_0}(\tau_B>ik_0) \leq (1/2)^i. $$
And then, 
\begin{eqnarray*}
\E_{x_0}(\tau_B) &\leq& \sum_{i\geq 0}  (i+1)k_0 \ \P_{x_0}  \left(\tau_B\in [ik_0;(i+1 )k_0 [ \right) \\
&\leq& \sum_{i\geq 0}  (i+1)k_0 \ \P_{x_0}  \left(\tau_B > (i+1 )k_0  \right) \\
& \leq&   c' k_0 \\
&\leq& c'' \nu(B)^{2/d}.\\
\end{eqnarray*}
Finally, we well retrieve  the second inequality (\ref{tempsperco}) of  proposition \ref{tau sur perco}.
\\

{\it{Acknowledgments}}: The authors would like to thank Noam Berger and Pierre Mathieu for their  comments on earlier version of the paper.

\nocite{*}
\bibliographystyle{plain}
\bibliography{biblioancree}

Thierry Delmotte
\\ \address{ Universit\'e Paul Sabatier  \\ 
Institut de Math\'ematiques de Toulouse  
\\ route de Narbonne
\\ 31400 Toulouse }
\\  \email{Thierry.Delmotte@math.ups-tlse.fr }
\\ \urladdr{http://www.math.univ-toulouse.fr/~delmotte}

Cl\'ement RAU
\\ \address{ Universit\'e Paul Sabatier  \\
Institut de Math\'ematiques de Toulouse
\\ route de Narbonne
\\ 31400 Toulouse }
\\  \email{rau@math.ups-tlse.fr}
\\ \urladdr{http://www.math.univ-toulouse.fr/~rau/}

\end{document}